\documentclass[preprint,12pt,authoryear]{elsarticle}

\usepackage{amssymb}
\usepackage{makecell} 
\usepackage{geometry,xcolor,tabularx,nccmath}
\usepackage{booktabs,caption,subcaption,mathtools,multirow,graphicx}

\usepackage{amsmath,amsthm,enumitem}

\usepackage[ruled,vlined,linesnumbered]{algorithm2e}

\usepackage{lineno}
\usepackage{tikz}
\usepackage{moreverb,url}
\usepackage{rotating}
\usepackage{pgfplots}

\usepackage[colorlinks=true,bookmarksopen=true,bookmarksnumbered=true,citecolor=red,urlcolor=red]{hyperref}
\usepackage[nameinlink,capitalise]{cleveref}
\usepackage{threeparttable}

\usepackage{tikz}
\usepackage{pgfplots}
\pgfplotsset{compat=1.18}
\usetikzlibrary{shapes.geometric,arrows.meta,positioning,fit,calc}

\definecolor{mplblue}{HTML}{1F77B4}    % BBB
\definecolor{mplorange}{HTML}{FF7F0E}  % CTA
\definecolor{mplgreen}{HTML}{2CA02C}   % PACE
\definecolor{mplred}{HTML}{D62728}     % KAT
\definecolor{mplpurple}{HTML}{9467BD}  % MARTA
\definecolor{mplbrown}{HTML}{8C564B}   % MTA

\journal{Transportation Research Part C: Emerging Technologies}

\begin{document}

\begin{frontmatter}

\title{Bus Fleet Electrification Under Capital Cost and Scheduling Constraints: A Five-Agency Case Study}

\author[1,2]{Hadi Bhidya}
\ead{hbhidya@anl.gov}

\author[1,3]{Taner Cokyasar\corref{cor1}}
\ead{tcokyasar@anl.gov}
\ead{tcokyasar@tamu.edu}

\author[1]{Omer Verbas}
\ead{omer@anl.gov}

\cortext[cor1]{Corresponding author}

\address[1]{Argonne National Laboratory, 9700 S. Cass Ave, Lemont, IL, 60439, USA}
\address[2]{The University of Tennessee, Knoxville, 1331 Circle Park Dr, Knoxville, TN, 37916, USA}
\address[3]{Texas A\&M University, Fermier Hall, 106 Ross St., College Station, TX, 77843, USA}

\begin{abstract}
As transit agencies consider bus fleet electrification, understanding the efficiency and cost of replacing diesel buses (DBs) with battery electric buses (BEBs) is critical. To evaluate this, this study applies a mixed-fleet optimization model—integrating scheduling, charging, and fleet composition decisions—across five agencies: Santa Monica's Big Blue Bus (BBB), the Chicago Transit Authority (CTA), Knoxville Area Transit (KAT), the Metropolitan Atlanta Rapid Transit Authority (MARTA), and Manhattan's Metropolitan Transportation Authority (MTA) bus service. By calculating electric fleet share, the BEB/DB replacement ratio, transit-link density, and vehicle activity-time allocation, the study finds that while optimized fleets remain majority-electric, vehicle substitution is rarely one-to-one. Average replacement ratios range from 1.101 for CTA to 1.245 for KAT, with higher transit-link density networks like CTA and MTA requiring fewer replacement buses per diesel bus displaced than lower-density networks like MARTA and KAT. While these relationships are descriptive rather than causal, non-revenue vehicle activity may help explain the differences. By shifting the focus from simple electric fleet share to diesel replacement efficiency, this multi-agency comparison demonstrates that transit agencies should use the replacement ratio to accurately forecast additional fleet capacity requirements and avoid the costly assumption of strict one-to-one vehicle substitution.
\end{abstract}


\end{frontmatter}

% =============================================================================
\section{Introduction}\label{sec:intro}
% =============================================================================

Transit agencies are increasingly evaluating battery electric buses (BEBs) as part of long-term fleet replacement planning \cite{perumal_2022,he_2023_transition}. BEBs may offer reductions in fuel and maintenance expenditures, but they differ from diesel buses (DBs) in both vehicle cost and daily operating requirements \cite{johnson_2020,perumal_2022}. In particular, BEBs must be scheduled around usable battery range, charging duration, charger availability, charging location, and garage access \cite{li_2014,adler_mirchandani_2017,perumal_2022,bazarnovi_2026}. The relevant planning question is therefore not only whether BEBs can cover scheduled service, but how efficiently they can replace the DB fleet required to provide that service.

Battery and charging constraints can change the set of trips that a single vehicle can cover. A sequence of trips that is feasible for a DB may require charging, additional deadheading, schedule restructuring, or an additional vehicle when assigned to a BEB \cite{li_2014,liu_ceder_2020,cokyasar_2023_large_scale,davatgari_2024_ndo}. These constraints make fleet transition a joint fleet-composition, vehicle-scheduling, and charging problem rather than a direct substitution of one vehicle type for another \cite{rinaldi_2020,alvo_2021,he_2023b,bazarnovi_2026}.

This study evaluates that issue using the BEB/DB replacement ratio. The ratio is defined as the number of BEBs selected in an optimized mixed fleet divided by the number of DBs displaced relative to an agency-specific diesel-only baseline. A value of 1.0 represents one BEB per DB displaced, while a value above 1.0 indicates that additional BEB fleet capacity is selected. Values below 1.0 are theoretically possible if the optimized mixed fleet requires fewer total vehicles than the diesel-only baseline, although such outcomes are not expected under the scheduling and charging constraints considered here. The metric complements electric fleet share, which describes fleet composition but does not directly indicate how efficiently BEBs replace DBs.

The analysis uses the integrated mixed BEB/DB scheduling and flexible-charging framework developed by \citet{bazarnovi_2026}. The model endogenously determines BEB and DB fleet counts while assigning scheduled trips and enforcing battery, charging, and charger-capacity constraints. The present study does not modify the formulation or solution procedure. Instead, it uses the resulting fleet and activity outputs to evaluate replacement efficiency across agencies and relative vehicle purchase-cost conditions.

Five transit systems are considered: Santa Monica's Big Blue Bus (BBB), the Chicago Transit Authority (CTA), Knoxville Area Transit (KAT), Metropolitan Atlanta Rapid Transit Authority (MARTA), and Manhattan's Metropolitan Transportation Authority (MTA). The DB purchase cost parameter is held fixed while the BEB purchase-cost parameter is varied to produce BEB/DB cost ratios from 1.0 to 1.4. Across these scenarios, the study evaluates electric fleet share, BEB/DB replacement ratio, transit-link density, and vehicle activity-time allocation. The spatial metric describes the concentration of each modeled network within a stop-based network footprint, while the activity metrics provide post-optimization diagnostics of revenue service, charging, deadheading, layover, waiting, and garage movements.

The contributions of this paper are threefold. First, it applies an established integrated mixed-fleet scheduling and charging model across five transit systems and a controlled range of relative BEB purchase costs. Second, it evaluates the BEB/DB replacement ratio as a system-level measure that distinguishes electric fleet share from diesel replacement efficiency. Third, it compares replacement outcomes with agency-level network-density and activity-time characteristics to investigate why identical vehicle-cost assumptions can produce different optimized fleet requirements across transit systems.

% =============================================================================
\section{Background and Literature Review}\label{sec:background}
% =============================================================================

Transit service planning is commonly represented as a sequence of interrelated decisions involving network design, frequency setting, timetabling, vehicle scheduling, and crew scheduling \cite{ceder_wilson_1986,guihaire_2008}. The vehicle scheduling problem (VSP) occurs after routes, frequencies, and timetabled trips have been established and determines how those trips are assigned to vehicles. In the conventional VSP, every scheduled trip must be covered, consecutive trips assigned to the same vehicle must be spatially and temporally compatible, and vehicle duties generally begin and end at an assigned garage. Common objectives include minimizing the required fleet, non-revenue travel, and operating cost \cite{freling_2001,bunte_2010}. The present study assumes that routes, frequencies, timetables, and garage assignments are fixed and focuses on the fleet outcomes produced by the vehicle-scheduling stage.

The introduction of limited-energy vehicles extends the conventional vehicle scheduling problem by adding fuel or battery-resource constraints. \citet{adler_mirchandani_2017} show how limited fuel capacity and fixed refueling locations can alter otherwise feasible vehicle duties. For BEBs, battery state of charge must remain feasible throughout each duty, and energy is consumed during both revenue and non-revenue movements. Early electric-bus studies incorporated limited battery energy, fast charging, battery exchange, and full or partial charging \cite{li_2014,zhu_2013,wen_2016}. These constraints mean that a trip sequence feasible for a DB may require charging, additional deadheading, schedule restructuring, or additional vehicles when assigned to a BEB.

Electric-bus planning therefore requires coordination among vehicle scheduling, charger location, and charging decisions \cite{perumal_2022}. Vehicle schedules determine when and where charging opportunities occur, while charger availability and charging duration affect whether a BEB can continue to its next trip or whether the duty must be split across multiple vehicles. Studies have optimized charging assignments, timing, duration, and charger use under state-of-charge, electricity-price, and station-capacity constraints \cite{abdelwahed_2020,bao_2023}. Integrated formulations have further combined vehicle scheduling with flexible charging, timetable shifting, charger deployment, and charging management \cite{duan_2023,he_2023b,gairola_2023,davatgari_2024_ebscl}. Together, this literature shows that BEB fleet requirements depend on the interaction between vehicle assignments and available charging opportunities.

Charging flexibility is especially important because it affects the amount of scheduled service that each BEB can cover. Opportunity and top-up charging can use scheduled dwell or layover periods and may reduce the need for vehicles to return to garages solely to restore battery energy \cite{hu_2021,hu_2022}. Partial charging can also permit a bus to obtain only the energy required for the next portion of its schedule rather than occupying a charger until the battery is full \cite{li_2020,bao_2023}. Conversely, limited charger counts or station capacities can create competition among buses, require waiting, or make otherwise feasible trip chains infeasible. \citet{liu_ceder_2020} explicitly minimize both the electric fleet size and the number of stationary chargers, demonstrating that these resources are interdependent. \citet{jovanovic_2021} similarly show that charging power and battery range affect required fleet size. These findings indicate that flexible and well-located charging can improve BEB productivity, whereas restrictive charging configurations can increase charging-related travel, waiting, and fleet requirements.

Related studies have integrated electric vehicle scheduling with crew scheduling, timetable design, multiple vehicle types, battery decisions, and long-term infrastructure deployment \cite{perumal_2021,xu_2023_timetable,yao_2020,alwesabi_2020,he_2023_transition}. These extensions demonstrate the interdependence of BEB planning decisions, although crew scheduling, timetable redesign, battery sizing, and long-term deployment remain outside the scope of the present fixed-service analysis.

Many electric-bus scheduling studies assume a fully electric fleet or fix the number of BEBs in advance. Mixed-fleet models instead allow BEBs and conventional buses to be assigned jointly. \citet{rinaldi_2020} formulate a mixed-fleet scheduling problem with electric and conventional or hybrid buses, while \citet{alvo_2021} consider a mixed electric--diesel fleet with restricted charger capacity and flexible charging. These models allow the level of electrification to emerge from vehicle costs, charging constraints, and operating requirements.

Mixed-fleet composition is also closely connected to capital and operating costs. BEBs may reduce fuel and maintenance expenditures, but their purchase prices remain influential inputs in fleet-investment analysis. \citet{johnson_2020} identify BEB and conventional-vehicle purchase prices as two of the most influential financial parameters, while mixed-fleet scheduling studies show that vehicle costs must be considered together with fleet size and vehicle productivity \cite{rinaldi_2020,alvo_2021}. Prior studies also demonstrate that operationally feasible electrification may require additional vehicles. \citet{cokyasar_2023_large_scale} report approximately 1.6 BEBs per displaced diesel bus under a tested 150-mile-range condition, while \citet{davatgari_2024_ndo} and \citet{bazarnovi_2026} show that range, next-day operability, garage-only charging, and non-revenue activity can increase fleet requirements.

Existing studies commonly report fleet size, electric fleet share, costs, charging resources, and non-revenue activity \cite{liu_ceder_2020,rinaldi_2020,davatgari_2024_ndo,bazarnovi_2026}, but they do not typically express fleet outcomes as the number of BEBs selected per DB displaced relative to a diesel-only baseline. This study uses the BEB/DB replacement ratio to make that relationship explicit and compares the resulting values across agencies. Transit-link density and vehicle activity-time allocation are then used as descriptive indicators of the spatial and operational conditions associated with differences in replacement efficiency.

Replacement efficiency may also differ across agencies even when vehicle costs and technology assumptions are held constant. Prior scheduling studies show that feasible BEB duties depend on trip compatibility, deadheading, garage and charger access, layover availability, and charging constraints \cite{adler_mirchandani_2017,liu_ceder_2020,alvo_2021,davatgari_2024_ndo,bazarnovi_2026}. These conditions vary across transit systems and can affect the share of each vehicle's schedule devoted to revenue service. Activity-level outputs can help identify these mechanisms by separating revenue service from charging, waiting, deadheading, layover, pull-out, and pull-in activities \cite{cokyasar_2023_large_scale,davatgari_2024_ndo,bazarnovi_2026}. The present study uses transit-link density only as a descriptive indicator of agency-scale spatial concentration; prior literature supports the underlying scheduling mechanisms but does not establish a direct causal relationship between network density and BEB replacement efficiency.

Limited work expresses optimized fleet outcomes as BEBs selected per DB displaced, compares that measure across agencies under common relative purchase-cost assumptions, and relates the results to spatial and activity-time characteristics. This study addresses that gap by applying the framework developed by \citet{bazarnovi_2026} to five transit systems, calculating the BEB/DB replacement ratio relative to agency-specific diesel-only baselines, and comparing replacement outcomes with transit-link density and vehicle activity-time allocation.

% =============================================================================
\section{Methodology}\label{sec:method}
% =============================================================================

The analysis applies the integrated mixed BEB/DB scheduling and flexible-charging model developed by \citet{bazarnovi_2026}. For each agency and purchase-cost scenario, the model jointly determines fleet composition, vehicle schedules, and charging decisions while producing activity-level vehicle records. These outputs are post-processed to calculate replacement, transit-link-density, and vehicle activity-time metrics.

\subsection{Model Selection and Analytical Scope}

The replacement-efficiency experiment requires a model in which BEB and DB fleet sizes emerge as decision outcomes rather than fixed inputs. The integrated framework developed by \citet{bazarnovi_2026} jointly determines mixed-fleet composition, vehicle schedules, and flexible charging decisions while ensuring coverage of a fixed set of timetabled trips. The model accounts for vehicle purchase and operating costs, battery and state-of-charge constraints, charging location and duration, partial charging, multiple charger types, and charger-capacity limits.

This structure is necessary because models that assume a fully electric fleet or prescribe the number of BEBs cannot represent the endogenous retention of DBs or changes in fleet composition across purchase-cost scenarios. Likewise, charging-scheduling and infrastructure-location models do not necessarily determine BEB and DB fleet requirements within the same daily scheduling problem \cite{abdelwahed_2020,liu_ceder_2020,alvo_2021,bao_2023,davatgari_2024,bazarnovi_2025_charging_location}. The \citet{bazarnovi_2026} framework therefore provides the decision structure needed to calculate the number of BEBs selected per DB displaced.

The present study does not modify the model formulation, objective function, charging representation, or solution algorithm. Instead, optimized fleet compositions, vehicle schedules, and activity records are post-processed to calculate replacement, spatial, and activity-time metrics. Unless otherwise specified, vehicle characteristics, battery and state-of-charge parameters, charging-power assumptions, energy-consumption rates, fuel and electricity prices, operating-cost parameters, and scheduling constraints are retained from the base-case configuration reported by \citet{bazarnovi_2026}. The agency service inputs and BEB purchase-cost parameter are varied, while the underlying vehicle, energy, and charging-performance assumptions remain unchanged.

\subsection{Case Study Systems}

The analysis is conducted for five transit systems: BBB, CTA, KAT, MARTA, and MTA. The systems provide variation in optimized diesel baseline fleet size, network footprint, transit-link density, and activity-time allocation. \Cref{tab:case_study_systems} summarizes the five systems and their diesel-only model
outputs.

\begin{table}[!ht]
\centering
\caption{Case study transit systems}
\label{tab:case_study_systems}
\begin{tabularx}{\textwidth}{l X r}
\hline
Region & Transit Agency & Modeled Diesel-Only Fleet \\
\hline
Santa Monica, CA & Big Blue Bus (BBB) & 126 \\
Chicago, IL & Chicago Transit Authority (CTA) & 1,543 \\
Knoxville, TN & Knoxville Area Transit (KAT) & 54 \\
Atlanta, GA & Metropolitan Atlanta Rapid Transit Authority (MARTA) & 426 \\
Manhattan & Metropolitan Transportation Authority (MTA), Manhattan bus service & 595 \\
\hline
\end{tabularx}
\end{table}

For each agency, the diesel-only baseline fleet is calculated as the number of unique DBs required by the optimization model to cover the fixed scheduled service. A Wednesday weekday schedule is used for each system. The same agency-specific diesel-only baseline is then used as the reference when calculating the number of DBs displaced in every mixed BEB/DB cost scenario.

Scheduled trips, stops, stop times, service calendars, and available route shapes were obtained from static General Transit Feed Specification (GTFS) feeds published by the respective transit agencies \cite{mobilitydata_gtfs_2026,bbb_gtfs_2026,cta_gtfs_2026,kat_gtfs_2026,marta_gtfs_2026,mta_manhattan_gtfs_2026}. GTFS provides a common structure for representing scheduled transit service across agencies. The BBB, KAT, MARTA, and MTA feeds were accessed on February 10, 2026; January 23, 2026; February 18, 2026; and March 31, 2026, respectively. The CTA analysis used a GTFS feed dated September 12, 2019. Garage locations were obtained from the input data used by the scheduling model and are displayed with the modeled transit networks.

\begin{figure}[!ht]
\centering

% Top row: 2 centered maps
\begin{subfigure}{0.31\textwidth}
    \centering
    \includegraphics[width=\textwidth,height=0.16\textheight,keepaspectratio]{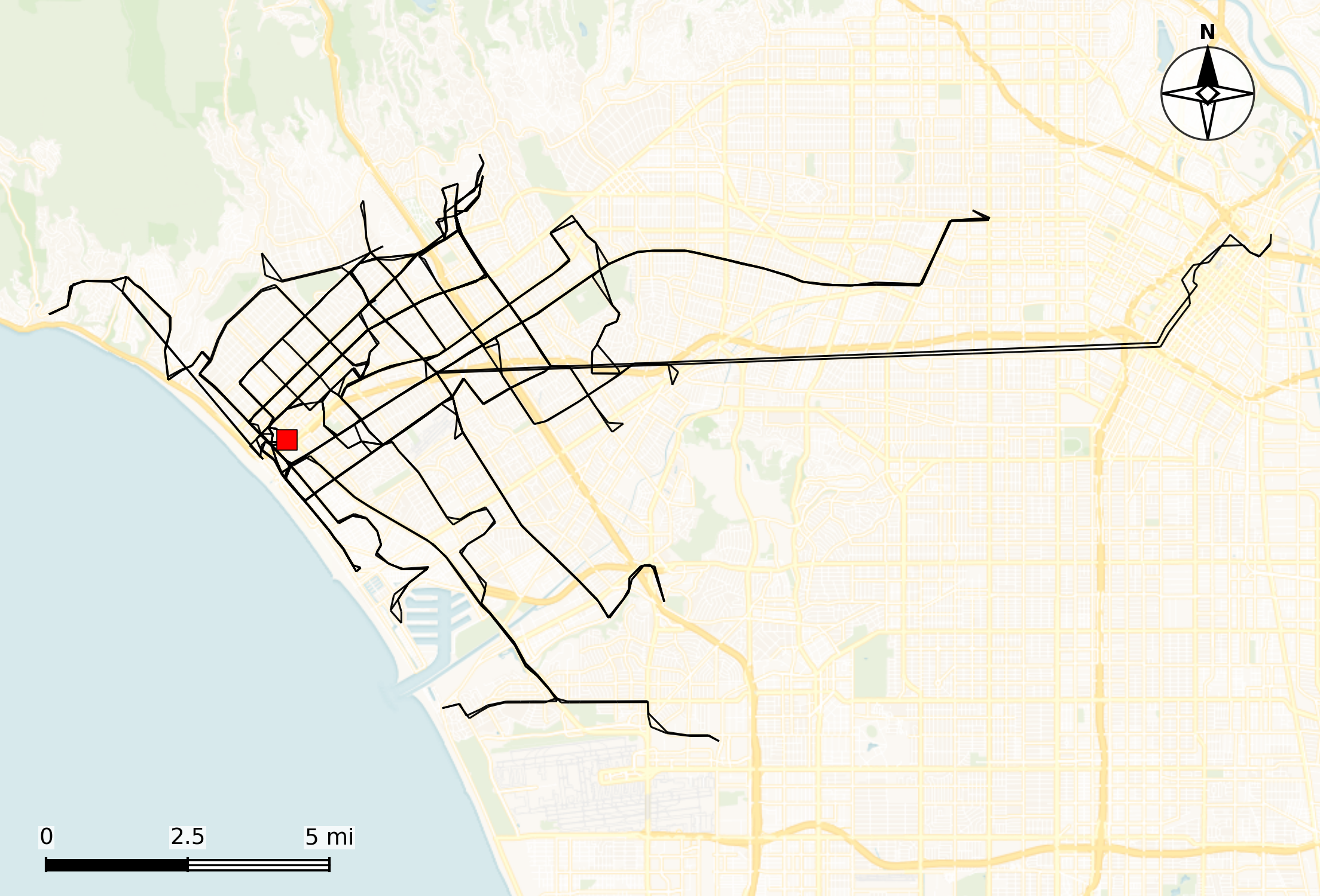}
    \caption{BBB}
    \label{fig:bbb_map}
\end{subfigure}
\hspace{0.04\textwidth}
\begin{subfigure}{0.31\textwidth}
    \centering
    \includegraphics[width=\textwidth,height=0.16\textheight,keepaspectratio]{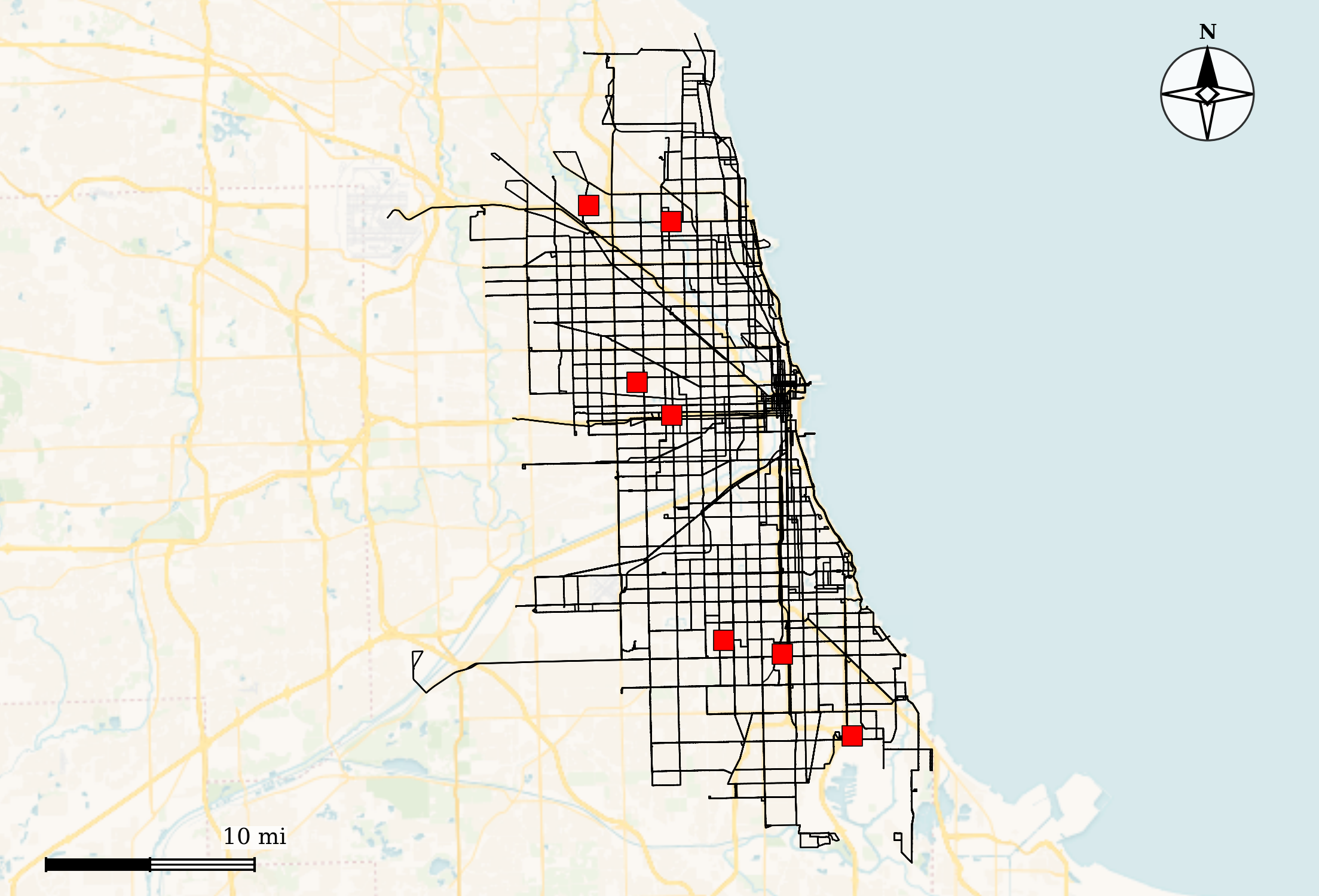}
    \caption{CTA}
    \label{fig:cta_map}
\end{subfigure}

% PACE removed
% \hfill
% \begin{subfigure}{0.31\textwidth}
%     \centering
%     \includegraphics[width=\textwidth,height=0.16\textheight,keepaspectratio]{transit_agency_maps/chicago-pace.png}
%     \caption{PACE}
%     \label{fig:pace_map}
% \end{subfigure}

\vspace{0.3em}

% Bottom row: 3 maps
\begin{subfigure}{0.31\textwidth}
    \centering
    \includegraphics[width=\textwidth,height=0.16\textheight,keepaspectratio]{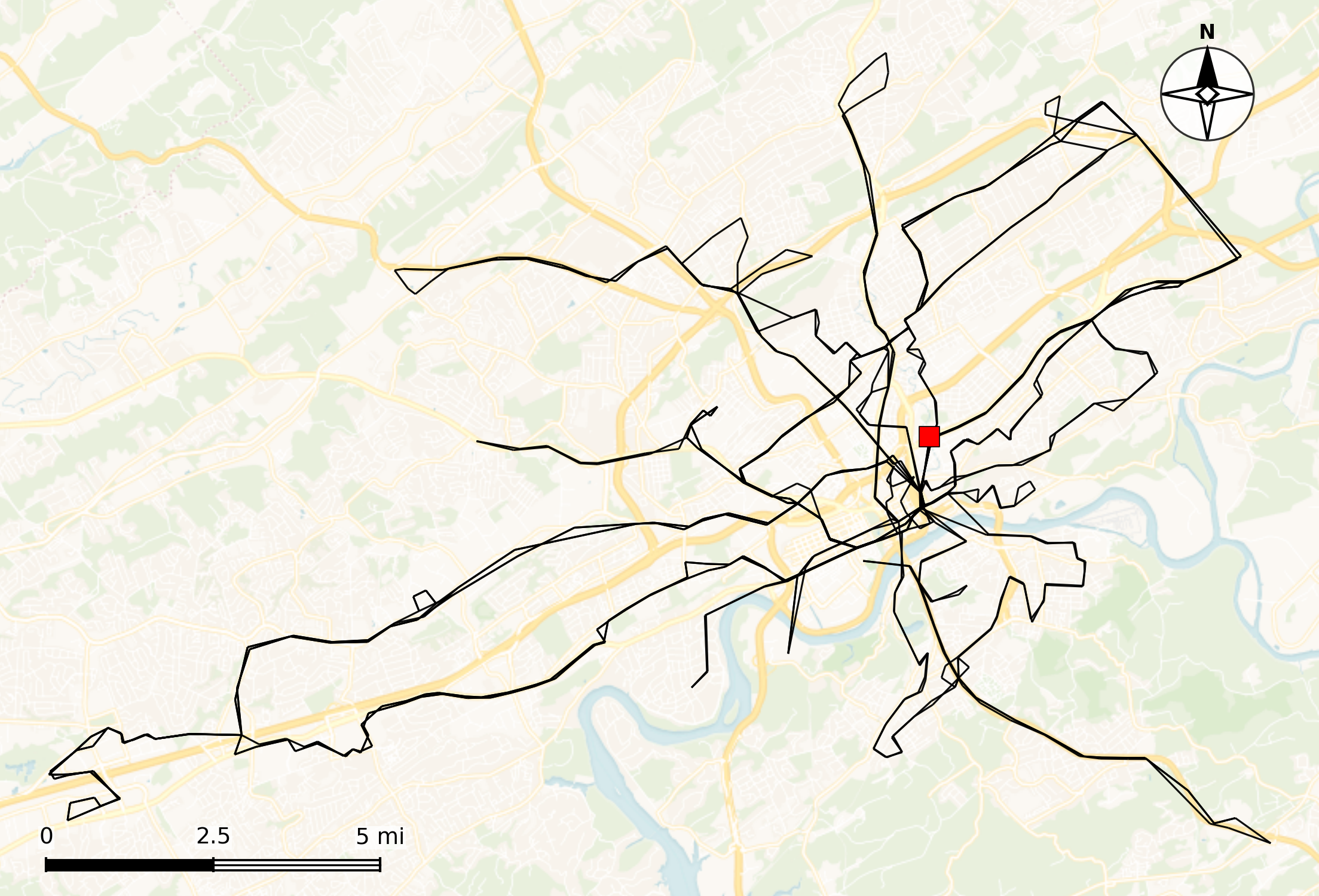}
    \caption{KAT}
    \label{fig:kat_map}
\end{subfigure}
\hfill
\begin{subfigure}{0.31\textwidth}
    \centering
    \includegraphics[width=\textwidth,height=0.16\textheight,keepaspectratio]{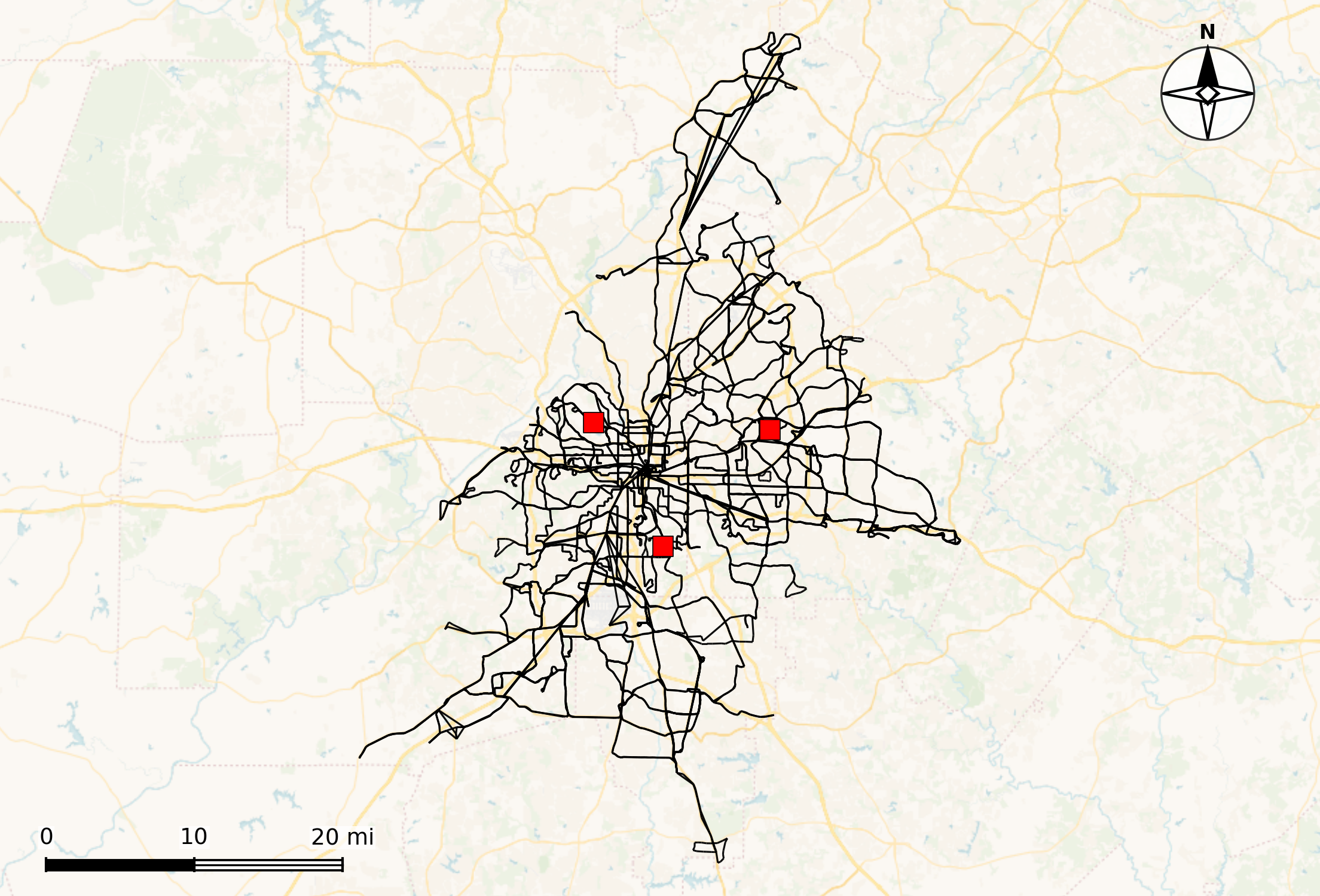}
    \caption{MARTA}
    \label{fig:marta_map}
\end{subfigure}
\hfill
\begin{subfigure}{0.31\textwidth}
    \centering
    \includegraphics[width=\textwidth,height=0.16\textheight,keepaspectratio]{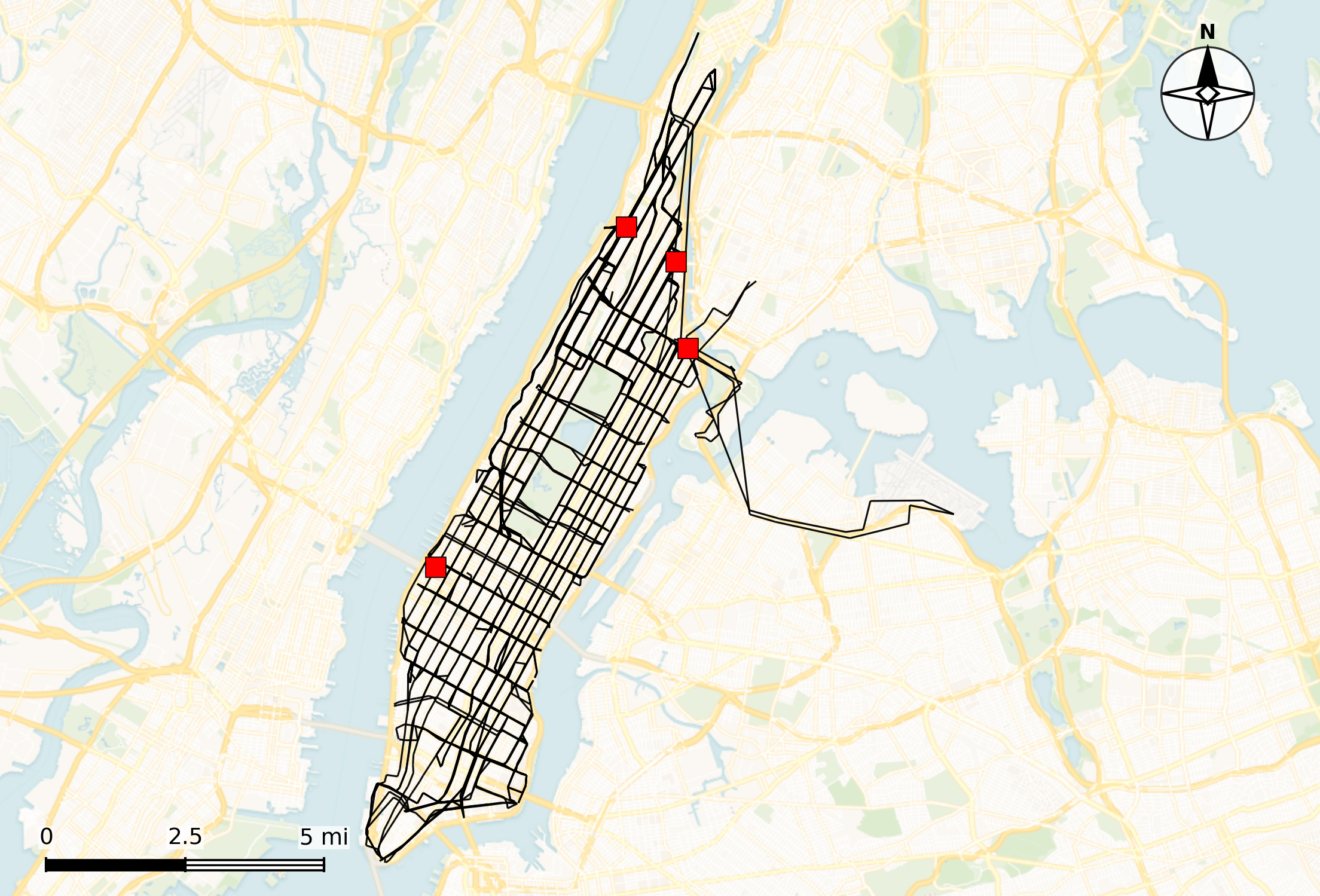}
    \caption{MTA}
    \label{fig:mta_map}
\end{subfigure}

\caption{Transit network maps for the five case study systems. Black lines represent transit links or GTFS route shapes, and red squares indicate garage locations.}
\label{fig:transit_network_maps}
\end{figure}

% \FloatBarrier

\subsection{Charging Infrastructure Assumptions}

Charging was restricted to the actual agency garage locations represented in the scheduling-model input data. For methodological simplicity and consistency across the five case-study systems, charger capacity at each garage was estimated from the amount of scheduled service associated with that garage. Trips were grouped by agency and garage, and garages serving more unique scheduled trips were assigned greater charging capacity.

The charger-allocation procedure uses \(C_g\) to denote the number of chargers assigned to garage \(g\), \(N_g\) to denote the number of unique scheduled trips associated with that garage, \(A_{\mathrm{UB}}\) to denote the model's upper-bound allocation parameter, \(m\) to denote the common charger-count multiplier, and \(u_g\) to denote a garage-specific scaling value between 0.01 and 0.03. The charger count is calculated as:

\begin{ceqn}
\begin{equation}
\label{eq:garage_chargers}
C_g =
\max\left\{
1,\,
\operatorname{round}
\left(
N_g A_{\mathrm{UB}} u_g m
\right)
\right\}.
\end{equation}
\end{ceqn}

The same allocation procedure and a fixed seed of 1 were applied to all agencies and purchase-cost scenarios. Consequently, the garage-specific scaling values and resulting charger inventories remained constant throughout the BEB/DB cost sweep. This approach preserves realistic charging locations while applying a consistent service-based charger-allocation framework across agencies, thereby reducing the influence of differing charging assumptions on the cross-agency replacement-efficiency comparison.

\subsection{Cost Scenarios}

BEB and DB purchase prices are influential and uncertain inputs in fleet-investment analysis. \citet{johnson_2020} identifies the purchase prices of the BEB and the displaced conventional vehicle as among the most influential parameters in BEB financial analysis. Accordingly, this study evaluates a controlled range of relative BEB/DB purchase costs rather than relying on a single cost relationship.

The model objective includes both capital and operating costs. Vehicle capital costs are based on purchase prices, annualized over a fixed vehicle lifetime, and converted to a daily equivalent for inclusion in the optimization objective. Operating costs include energy, fuel, and maintenance expenditures and are calculated over total vehicle operating activity, including both revenue service and non-revenue movements. Unless otherwise specified, the operating-cost assumptions are retained from the base-case configuration reported by \citet{bazarnovi_2026}.

All five agencies are evaluated under the same purchase-cost scenarios. The DB purchase-cost parameter is held fixed at \$500,000, while the BEB purchase-cost parameter is varied from \$500,000 to \$700,000 in \$20,000 increments. This produces BEB/DB purchase-cost ratios from 1.0 to 1.4. Only the BEB purchase-cost parameter is varied across the sensitivity analysis; the DB purchase cost, diesel fuel price, electricity price, maintenance costs, and all other operating-cost parameters remain fixed.

The electric fleet share is not prescribed for any purchase-cost scenario. Instead, the model is allowed to select the cost-minimizing fleet composition subject to an electrification lower bound of 1\% and an upper bound of 100\%. Within these bounds, the number of BEBs and DBs is determined endogenously based on the modeled capital costs, operating costs, scheduling requirements, battery constraints, and charging capacity. The electrification percentages reported in \Cref{fig:cost_ratio_electrification} are therefore optimization outcomes rather than scenario inputs.

The cost sweep is designed to isolate how changes in the relative capital cost of BEBs affect optimized fleet composition, total system cost, and replacement efficiency. Lower cost ratios represent cases in which the BEB purchase cost is closer to the DB purchase cost, while higher ratios represent a larger BEB purchase-cost premium. Results are therefore interpreted primarily in terms of the BEB/DB purchase-cost ratio rather than as forecasts based on any single absolute vehicle price.

\subsection{Fleet Replacement Metrics}

To interpret the optimized mixed-fleet outcomes, this study defines the BEB/DB replacement ratio for each agency and cost scenario. The ratio compares the number of BEBs selected in mixed scenario \(s\) with the number of DBs displaced relative to the diesel-only baseline. Here, \(N_{\mathrm{BEB},s}\) denotes the number of BEBs selected in mixed scenario \(s\), \(N_{\mathrm{DB},s}\) denotes the number of DBs remaining in that scenario, and \(N_{\mathrm{DB},0}\) denotes the number of DBs in the diesel-only baseline. The replacement ratio is calculated as:

\begin{ceqn}
\begin{equation}
\label{eq:replacement_ratio}
R_s =
\frac{N_{\mathrm{BEB},s}}
{N_{\mathrm{DB},0} - N_{\mathrm{DB},s}} .
\end{equation}
\end{ceqn}

The denominator in \Cref{eq:replacement_ratio} represents the number of DBs displaced relative to the diesel-only baseline, while the numerator represents the number of BEBs selected in the optimized mixed fleet. A value of 1.0 indicates one BEB per DB displaced, while values above 1.0 indicate that more BEBs are selected than DBs displaced.

Electric fleet share is also reported for each scenario as the percentage of optimized fleet vehicles that are BEBs. Electric fleet share is used to show the extent of BEB adoption in the optimized fleet, while the replacement ratio is used to measure how efficiently BEBs replace diesel service.

\subsection{Stop-Based Network Footprint and Transit-Link Density}

Transit-link density is calculated using the static GTFS stops and the transit network represented in the model. The GTFS stop coordinates are first used to construct a consistent spatial footprint for each agency.

To obtain a consistent spatial denominator across agencies, this study constructs a concave boundary around each agency's GTFS stop coordinates. Concave or alpha-shape boundaries generalize the convex hull by allowing the estimated polygon to follow nonconvex features in a point distribution more closely \cite{edelsbrunner_1983}. The resulting polygon is interpreted as a stop-based network footprint rather than a passenger walking catchment or official agency service boundary.

The alpha parameter was calibrated separately for each agency through visual inspection of candidate boundaries. The selected values produced coherent stop-based network footprints while limiting excessive inclusion of stop-free areas and fragmentation of peripheral service corridors. Because the resulting footprints depend on this calibration, transit-link density is interpreted as a descriptive measure rather than an exact estimate of official service-area density.

Transit-link density compares the total mileage of unique transit links in the model network with the area of the stop-based network footprint. Here, \(M_{\mathrm{link}}\) denotes the total mileage of unique transit links, and \(A_{\mathrm{footprint}}\) denotes the footprint area in square miles. Unique links are used to avoid repeatedly counting the same physical network segments across multiple scheduled trips. Transit-link density is calculated as:

\begin{ceqn}
\begin{equation}
\label{eq:network_density}
D =
\frac{M_{\mathrm{link}}}
{A_{\mathrm{footprint}}}.
\end{equation}
\end{ceqn}

Transit-link density represents the mileage of modeled transit links per square mile of the stop-based network footprint. This measure summarizes the spatial concentration of each modeled network, but it does not directly measure service frequency, passenger accessibility, trip-to-trip compatibility, charger accessibility, or deadheading. The relationship between this metric and replacement efficiency is therefore interpreted descriptively and examined alongside activity-level scheduling outputs.

\subsection{Vehicle Time Allocation Metrics}

Vehicle activity records produced by the optimized schedules are used to compare how BEBs and DBs allocate scheduled time across operating activities \cite{bazarnovi_2026}. Each activity record includes a vehicle identifier, bus type, activity classification, start time, and end time. Activity duration is calculated from the start and end times and converted to hours.

The reported activity categories are revenue trip, charging, deadhead to station, deadhead to next trip, deadhead from station, layover at station, layover at trip, waiting, pull-out, and pull-in. For each agency and bus type, activity durations are pooled across all buses and all mixed-fleet scenarios containing BEBs. The share assigned to each activity is calculated by dividing its pooled duration by the total pooled scheduled activity duration for that agency and bus type. Thus, each hour of scheduled vehicle activity receives equal weight, and vehicles or scenarios contributing more scheduled activity time contribute proportionally more to the resulting shares.

The activity shares are used as post-optimization diagnostics of vehicle productivity. Charging and waiting occupy portions of the schedule during which a BEB does not cover revenue trips. Deadheading and pull-out/pull-in movements also consume scheduled time and, for BEBs, battery energy \cite{cokyasar_2023_large_scale,davatgari_2024_ndo}. Layover is interpreted separately because it may represent non-revenue time but may also provide an opportunity to accommodate charging without requiring another vehicle. The activity categories are therefore not treated as equivalent burdens, and their relationship with replacement efficiency is interpreted in the context of the complete optimized schedule.

% =============================================================================
\section{Results}\label{sec:results}
% =============================================================================

\subsection{Fleet Electrification Across Cost Scenarios}

The first set of results evaluates how optimized fleet composition changes as the BEB/DB capital cost ratio increases. \Cref{fig:cost_ratio_electrification} shows, for each agency, the change in total cost relative to the case where BEB and DB purchase costs are equal, along with the optimized electrification rate.

% ------------------------------------------------------------
% Cost-ratio plot macro
% Usage:
% \CostRatioPlot{csv_path}{left_y_label}{right_y_label}
% ------------------------------------------------------------

\newcommand{\CostRatioPlot}[3]{%
\makebox[\linewidth][c]{%
\begin{tikzpicture}
\begin{axis}[
    width=\linewidth,
    height=0.58\linewidth,
    xmin=0.98,
    xmax=1.42,
    xlabel={BEB / DB Cost Ratio},
    ylabel={#2},
    axis y line*=left,
    axis x line*=bottom,
    ylabel style={font=\small\bfseries, text=blue!70!black},
    xlabel style={font=\small\bfseries},
    tick label style={font=\footnotesize},
    yticklabel style={text=blue!70!black},
    xtick={1.00,1.10,1.20,1.30,1.40},
    xticklabel style={
        /pgf/number format/fixed,
        /pgf/number format/precision=2
    },
    scaled ticks=false,
    grid=major,
    major grid style={black!10},
    axis line style={black},
    tick align=outside,
    enlargelimits=false,
    clip=false,
]
\addplot[
    blue!70!black,
    very thick,
    mark=o,
    mark size=2.4pt,
]
table[
    col sep=comma,
    x=ratio,
    y=cost_change
]{#1};

\addplot[
    blue!70!black,
    dotted,
    thick,
    mark=none,
]
coordinates {(0.98,0) (1.42,0)};
\end{axis}

\begin{axis}[
    width=\linewidth,
    height=0.58\linewidth,
    xmin=0.98,
    xmax=1.42,
    axis y line*=right,
    axis x line=none,
    ylabel={#3},
    ylabel style={font=\small\bfseries, text=red!75!black},
    tick label style={font=\footnotesize},
    yticklabel style={text=red!75!black},
    scaled ticks=false,
    grid=none,
    tick align=outside,
    enlargelimits=false,
    clip=false,
]
\addplot[
    red!75!black,
    very thick,
    dashed,
    mark=square*,
    mark size=2.3pt,
    nodes near coords={\pgfmathprintnumber[fixed,precision=0]{\pgfplotspointmeta}\%},
    point meta=explicit,
    every node near coord/.append style={
        font=\scriptsize,
        text=red!75!black,
        yshift=6pt
    },
]
table[
    col sep=comma,
    x=ratio,
    y=electrification,
    meta=electrification
]{#1};
\end{axis}
\end{tikzpicture}%
}%
}

\begin{figure}[!ht]
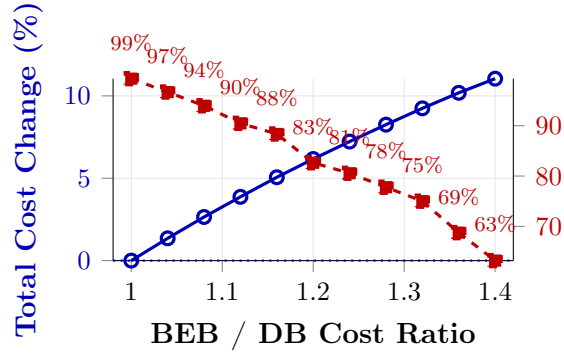
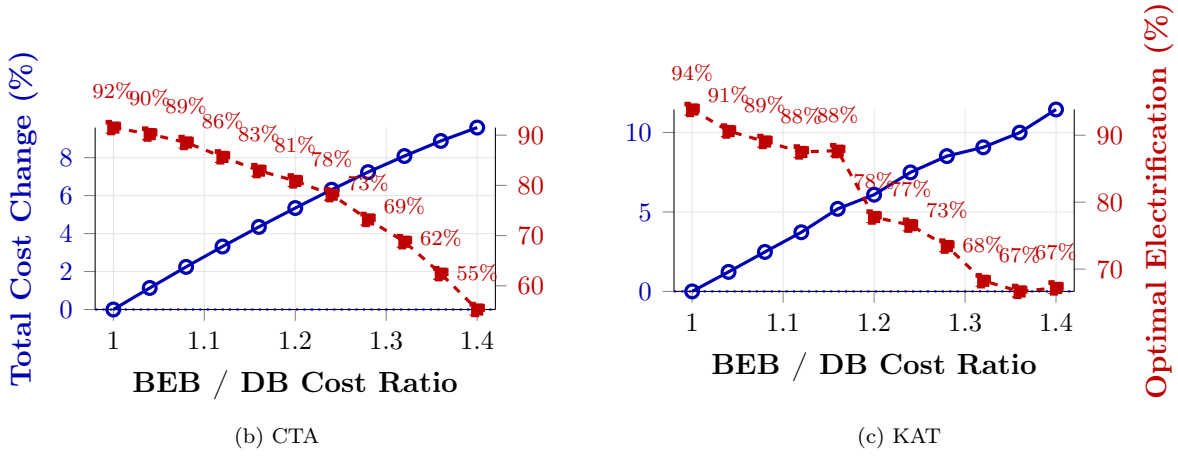
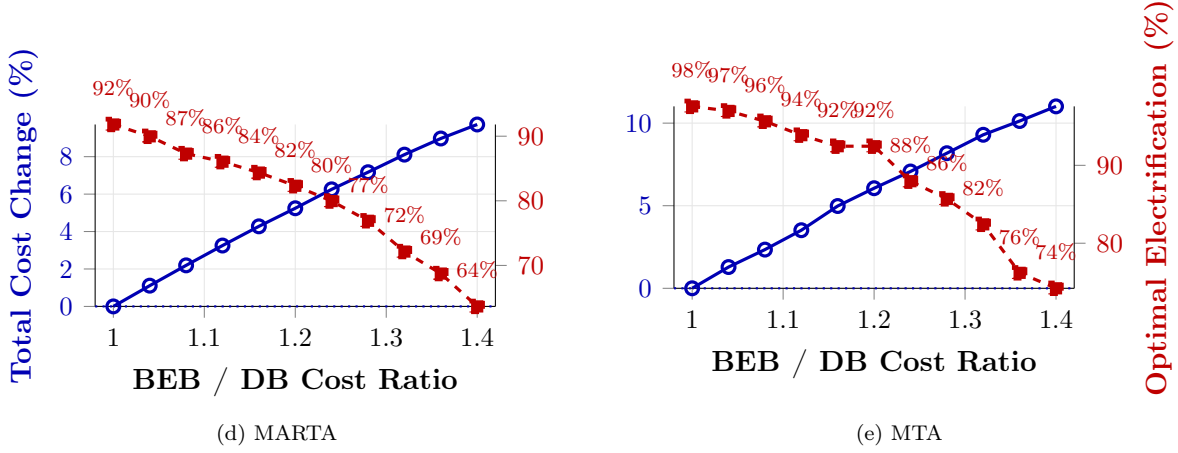

\centering

% Top row: 1 centered plot
\begin{subfigure}{0.455\textwidth}
    \centering
    \CostRatioPlot{cost_ratios/tikz_data/cost_ratio_bbb.csv}{Total Cost Change (\%)}{}
    \caption{BBB}
    \label{fig:cost_ratio_bbb}
\end{subfigure}

\vspace{0.4em}

% Middle row: 2 plots
\begin{subfigure}{0.455\textwidth}
    \centering
    \CostRatioPlot{cost_ratios/tikz_data/cost_ratio_cta.csv}{Total Cost Change (\%)}{}
    \caption{CTA}
    \label{fig:cost_ratio_cta}
\end{subfigure}
\hfill
\begin{subfigure}{0.455\textwidth}
    \centering
    \CostRatioPlot{cost_ratios/tikz_data/cost_ratio_kat.csv}{}{Optimal Electrification (\%)}
    \caption{KAT}
    \label{fig:cost_ratio_kat}
\end{subfigure}

\vspace{0.4em}

% PACE removed
% \begin{subfigure}{0.455\textwidth}
%     \centering
%     \CostRatioPlot{cost_ratios/tikz_data/cost_ratio_pace.csv}{Total Cost Change (\%)}{}
%     \caption{PACE}
%     \label{fig:cost_ratio_pace}
% \end{subfigure}

% Bottom row: 2 plots
\begin{subfigure}{0.455\textwidth}
    \centering
    \CostRatioPlot{cost_ratios/tikz_data/cost_ratio_marta.csv}{Total Cost Change (\%)}{}
    \caption{MARTA}
    \label{fig:cost_ratio_marta}
\end{subfigure}
\hfill
\begin{subfigure}{0.455\textwidth}
    \centering
    \CostRatioPlot{cost_ratios/tikz_data/cost_ratio_mta.csv}{}{Optimal Electrification (\%)}
    \caption{MTA}
    \label{fig:cost_ratio_mta}
\end{subfigure}

\caption{Total cost increase and optimized electrification rate across BEB/DB cost ratios. Blue lines correspond to total cost change, and red dashed lines correspond to optimized electrification rate.}
\label{fig:cost_ratio_electrification}
\end{figure}
% \FloatBarrier

Across all five analyzed agencies, the same general pattern appears: as the BEB/DB purchase-cost ratio increases, total optimized cost rises and the optimized electric fleet share declines. Nevertheless, the five optimized fleets remain majority-electric across the full tested cost range. This indicates that, under the modeled operating-cost and charging assumptions, BEBs continue to be selected even as their purchase-cost premium increases.

The scale and rate of decline differ across systems. BBB begins nearly fully electric when BEB and DB purchase costs are equal and remains majority-electric at the highest tested ratio. KAT, MARTA, and MTA show the same general cost-response pattern but retain different electric fleet shares across the sweep. CTA decreases from approximately 92\% electric at a cost ratio of 1.0 to approximately 55\% at 1.4, representing one of the larger declines among the five analyzed systems. Despite that decline, BEBs remain a majority of CTA's optimized fleet throughout the tested range.

These results demonstrate that electric fleet share describes the composition of the optimized fleet but does not reveal how many BEBs are selected for each DB removed. The next section therefore evaluates the BEB/DB replacement ratio.

\subsection{BEB/DB Replacement Ratios Across Agencies}

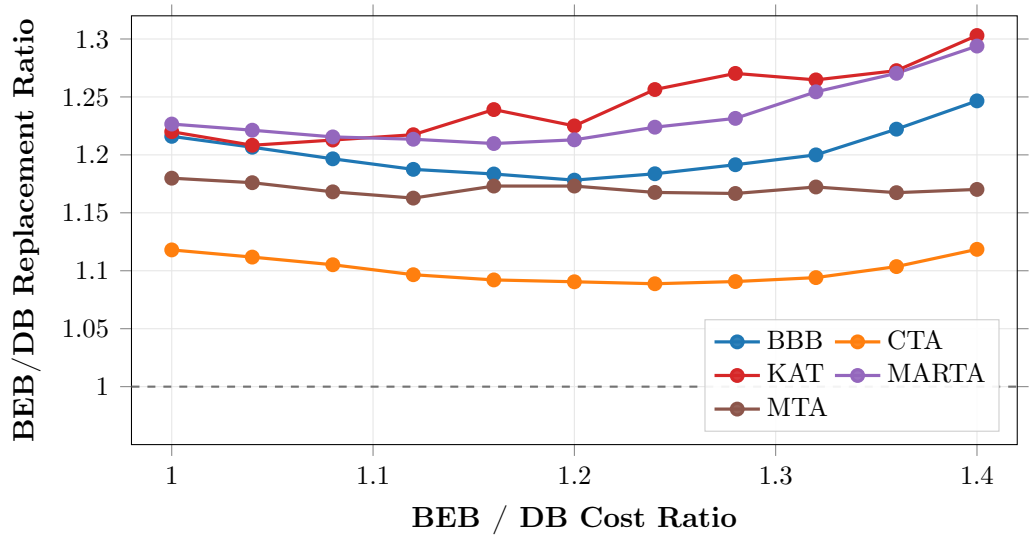
\begin{figure}[!ht]
\centering

\begin{tikzpicture}
\begin{axis}[
    width=0.88\textwidth,
    height=0.48\textwidth,
    xlabel={BEB / DB Cost Ratio},
    ylabel={BEB/DB Replacement Ratio},
    xlabel style={font=\small\bfseries},
    ylabel style={font=\small\bfseries},
    tick label style={font=\footnotesize},
    xmin=0.98,
    xmax=1.42,
    xtick={1.00,1.10,1.20,1.30,1.40},
    xticklabel style={
        /pgf/number format/fixed,
        /pgf/number format/precision=2
    },
    ymin=0.95,
    ymax=1.32,
    ytick={1.00,1.05,1.10,1.15,1.20,1.25,1.30},
    scaled ticks=false,
    grid=major,
    major grid style={black!10},
    axis line style={black},
    tick align=outside,
    legend style={
        at={(0.98,0.03)},
        anchor=south east,
        draw=black!20,
        fill=white,
        fill opacity=0.92,
        text opacity=1,
        font=\footnotesize
    },
    legend cell align={left},
    legend columns=2,
]

% One-to-one replacement reference line
\addplot[
    black!55,
    dashed,
    thick,
    mark=none,
    forget plot
]
coordinates {(0.98,1.0) (1.42,1.0)};

% BBB
\addplot[
    color=mplblue,
    very thick,
    mark=*,
    mark size=2.2pt,
]
table[
    col sep=comma,
    x=ratio,
    y=BBB
]{replacement_ratio/tikz_data/replacement_ratio_wide.csv};
\addlegendentry{BBB}

% CTA
\addplot[
    color=mplorange,
    very thick,
    mark=*,
    mark size=2.2pt,
]
table[
    col sep=comma,
    x=ratio,
    y=CTA
]{replacement_ratio/tikz_data/replacement_ratio_wide.csv};
\addlegendentry{CTA}

% PACE removed
% \addplot[
%     color=mplgreen,
%     very thick,
%     mark=*,
%     mark size=2.2pt,
% ]
% table[
%     col sep=comma,
%     x=ratio,
%     y=PACE
% ]{replacement_ratio/tikz_data/replacement_ratio_wide.csv};
% \addlegendentry{PACE}

% KAT
\addplot[
    color=mplred,
    very thick,
    mark=*,
    mark size=2.2pt,
]
table[
    col sep=comma,
    x=ratio,
    y=KAT
]{replacement_ratio/tikz_data/replacement_ratio_wide.csv};
\addlegendentry{KAT}

% MARTA
\addplot[
    color=mplpurple,
    very thick,
    mark=*,
    mark size=2.2pt,
]
table[
    col sep=comma,
    x=ratio,
    y=MARTA
]{replacement_ratio/tikz_data/replacement_ratio_wide.csv};
\addlegendentry{MARTA}

% MTA
\addplot[
    color=mplbrown,
    very thick,
    mark=*,
    mark size=2.2pt,
]
table[
    col sep=comma,
    x=ratio,
    y=MTA
]{replacement_ratio/tikz_data/replacement_ratio_wide.csv};
\addlegendentry{MTA}

\end{axis}
\end{tikzpicture}

\caption{BEB/DB replacement ratio across BEB/DB cost ratios by agency. The dashed horizontal line indicates one-to-one replacement; values above 1.0 indicate more BEBs than diesel buses displaced under the optimized fleet outcome.}
\label{fig:replacement_ratio_cost}
\end{figure}

\Cref{fig:replacement_ratio_cost} shows the BEB/DB replacement ratio for each agency across the tested BEB/DB cost ratios. This is the main replacement-efficiency result of the paper.

% \FloatBarrier

All reported mixed-fleet scenarios displaced at least one DB relative to the corresponding diesel-only baseline; therefore, the denominator was positive for every scenario included in the analysis.
The BEB/DB replacement ratio differs systematically across the five analyzed agencies. CTA maintains the lowest average replacement ratio and remains relatively stable across the
purchase-cost sweep, ranging approximately from 1.09 to 1.12. MTA has the second-lowest average ratio and also remains nearly flat, at approximately 1.16--1.18. These outcomes indicate that both systems maintain relatively stable replacement efficiency even as the BEB purchase-cost premium increases.

MARTA and KAT have the highest replacement ratios. Both agencies start above MTA and BBB and increase more noticeably as the BEB/DB cost ratio rises. Their replacement-ratio curves follow a similar pattern, suggesting that the optimized BEB fleets in these systems require more additional BEB capacity per DB displaced.

BBB occupies an intermediate position. Its replacement ratio initially decreases slightly and then rises at the higher cost ratios, reaching approximately 1.25 at the upper end of the sweep.

Average replacement ratios across all cost scenarios are 1.101 for CTA, 1.171 for MTA, 1.201 for BBB, 1.234 for MARTA, and 1.245 for KAT. Because all five agencies are evaluated under the same vehicle purchase-cost scenarios, the cross-agency variation cannot be attributed solely to the cost-ratio input. It indicates that differences in the scheduled service, trip connections, charging opportunities, garage relationships, and resulting vehicle assignments may contribute to the replacement outcomes.

\subsection{Transit-Link Density and Replacement Efficiency}

\begin{figure}[!ht]
\centering

% Top row: 2 centered maps
\begin{subfigure}{0.31\textwidth}
    \centering
    \includegraphics[width=\textwidth,height=0.16\textheight,keepaspectratio]{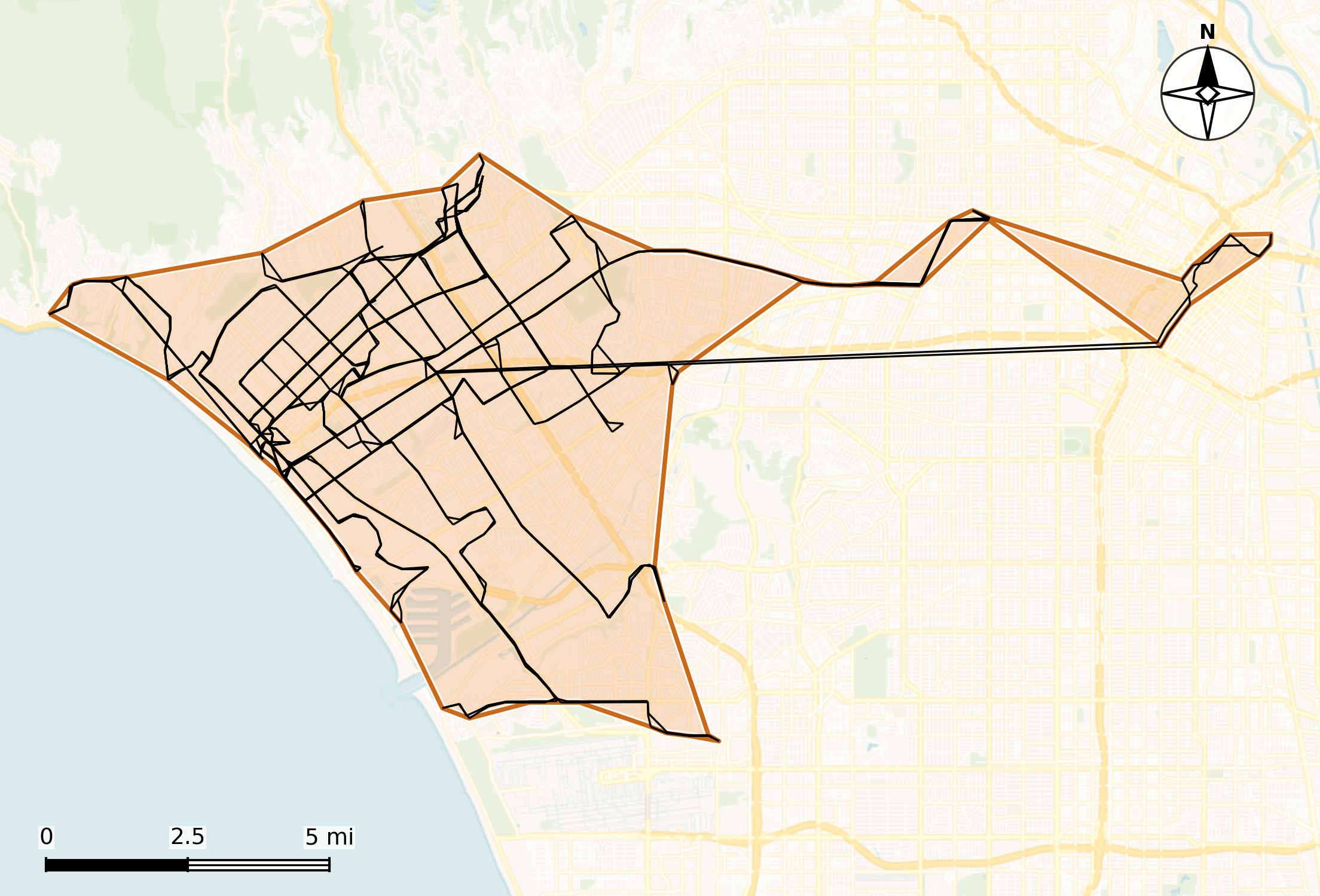}
    \caption{BBB}
    \label{fig:bbb_density_map}
\end{subfigure}
\hspace{0.04\textwidth}
\begin{subfigure}{0.31\textwidth}
    \centering
    \includegraphics[width=\textwidth,height=0.16\textheight,keepaspectratio]{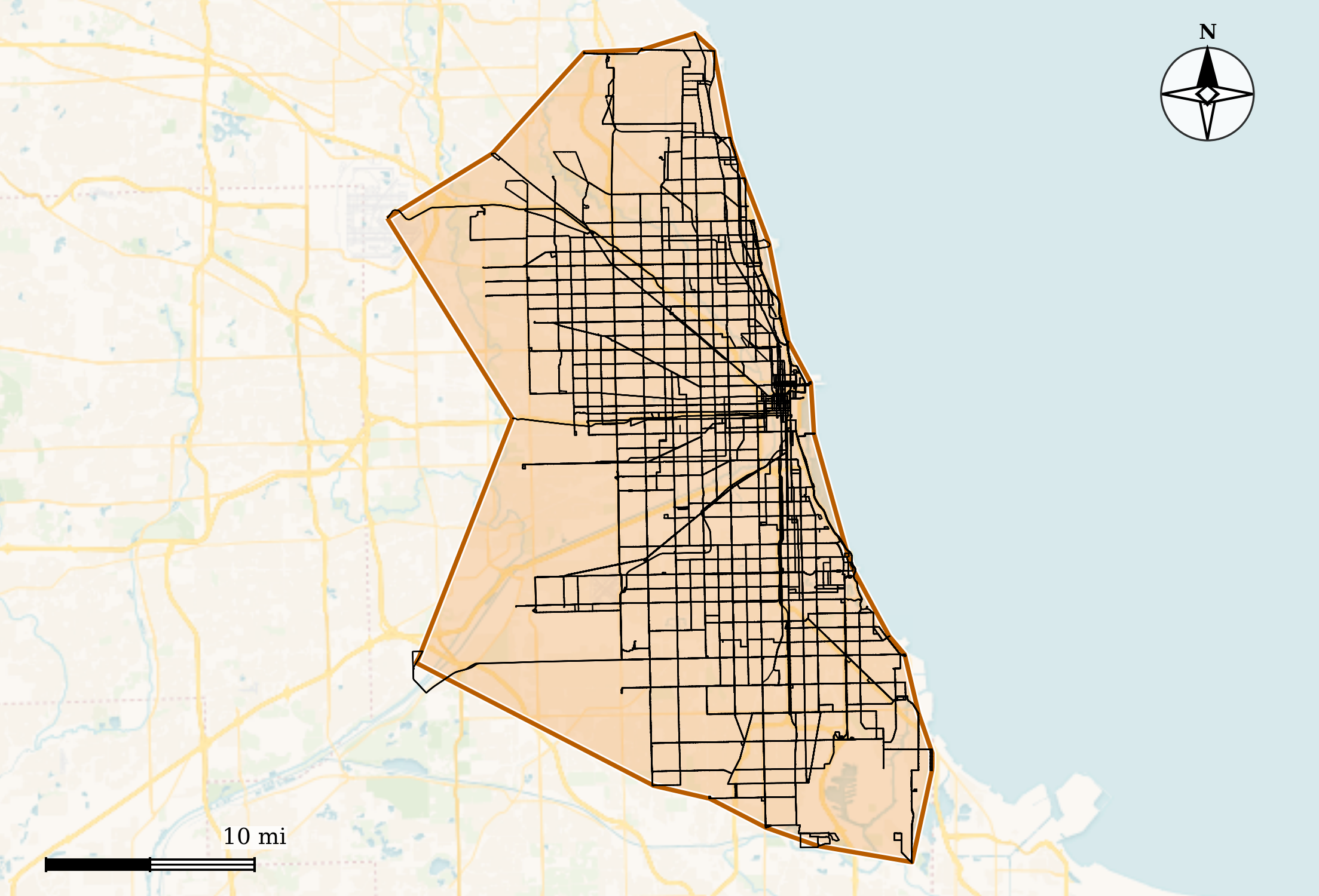}
    \caption{CTA}
    \label{fig:cta_density_map}
\end{subfigure}

% PACE removed
% \hfill
% \begin{subfigure}{0.31\textwidth}
%     \centering
%     \includegraphics[width=\textwidth,height=0.16\textheight,keepaspectratio]{density/chicago-pace_density_map.png}
%     \caption{PACE}
%     \label{fig:pace_density_map}
% \end{subfigure}

\vspace{0.3em}

% Bottom row: 3 maps
\begin{subfigure}{0.31\textwidth}
    \centering
    \includegraphics[width=\textwidth,height=0.16\textheight,keepaspectratio]{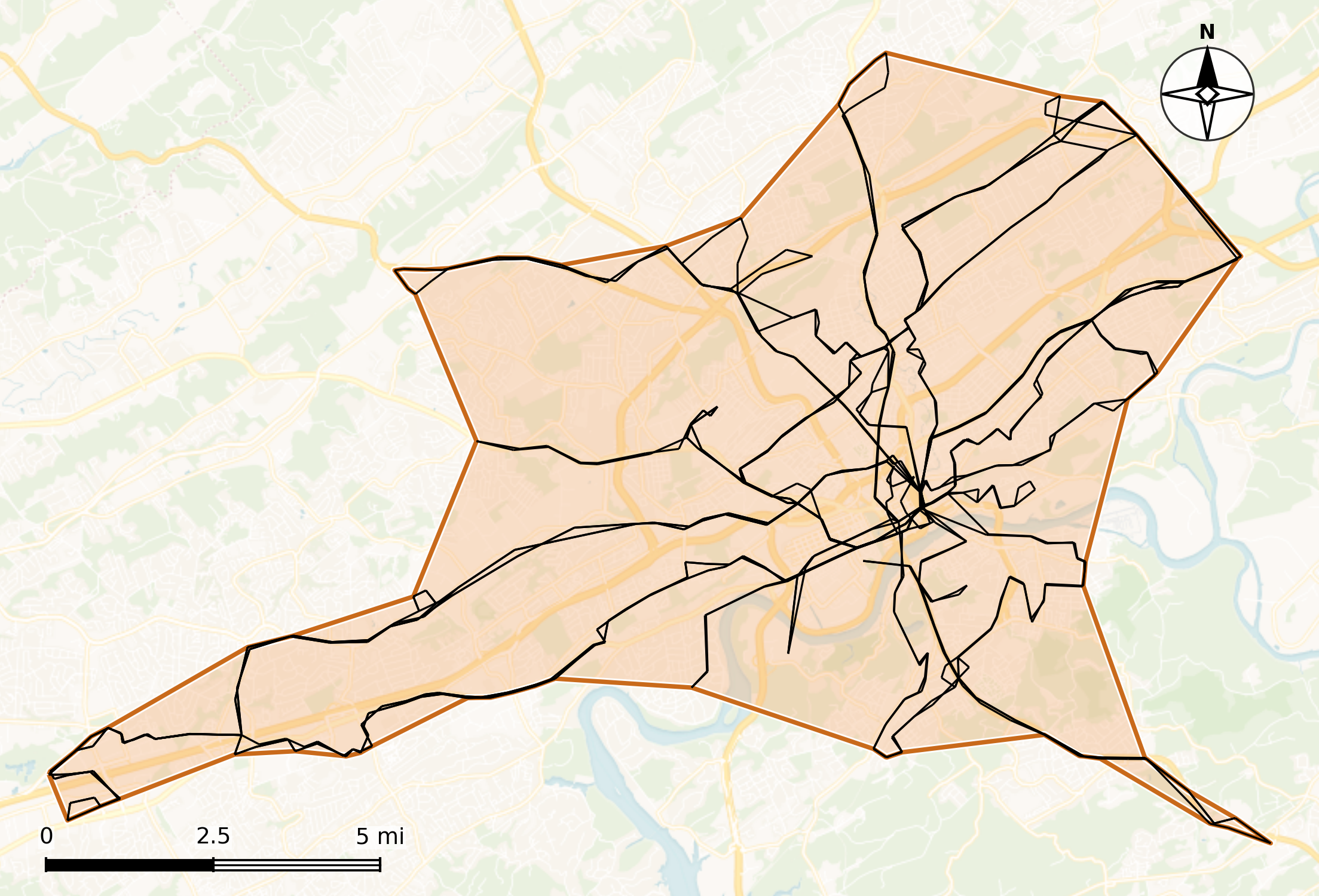}
    \caption{KAT}
    \label{fig:kat_density_map}
\end{subfigure}
\hfill
\begin{subfigure}{0.31\textwidth}
    \centering
    \includegraphics[width=\textwidth,height=0.16\textheight,keepaspectratio]{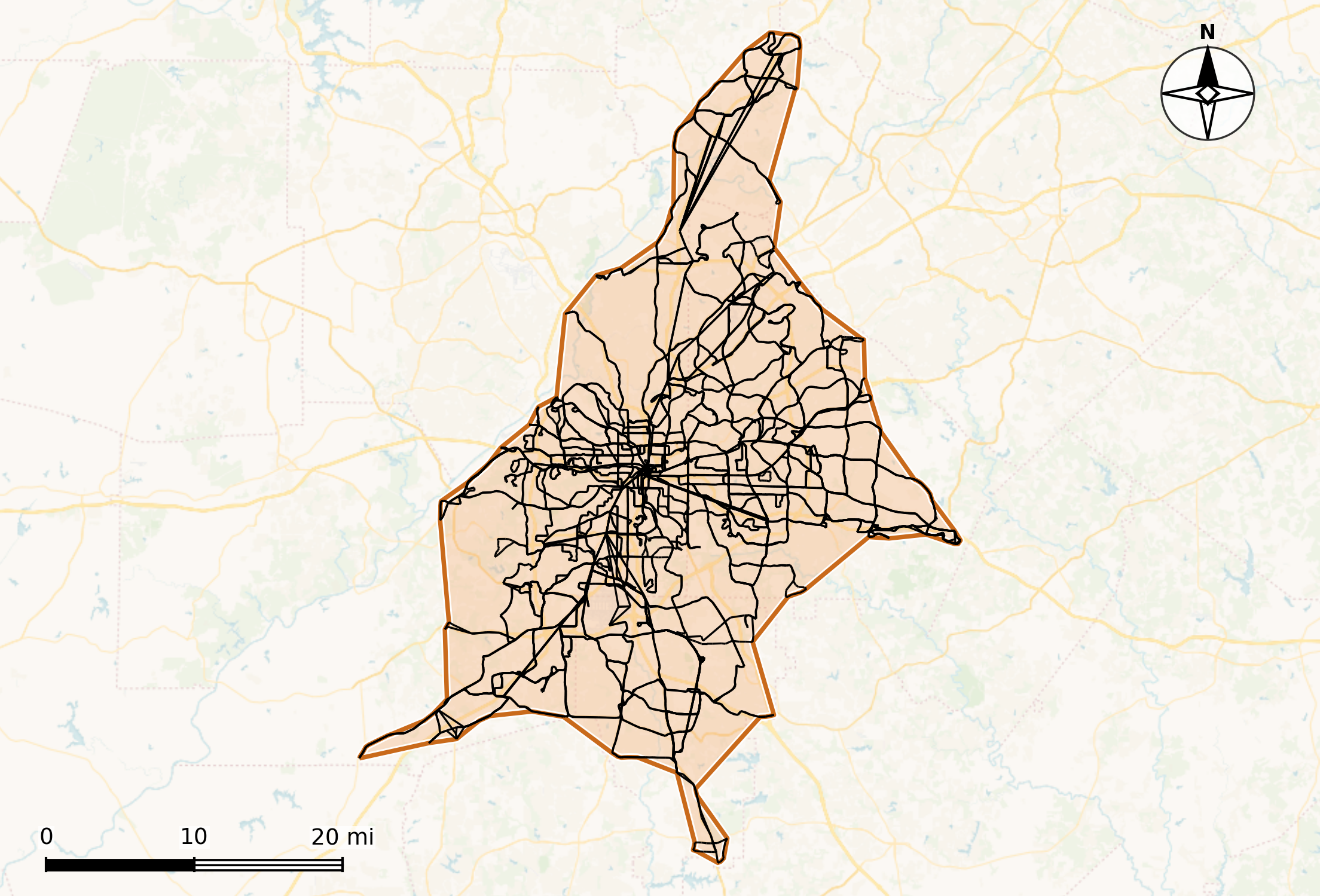}
    \caption{MARTA}
    \label{fig:marta_density_map}
\end{subfigure}
\hfill
\begin{subfigure}{0.31\textwidth}
    \centering
    \includegraphics[width=\textwidth,height=0.16\textheight,keepaspectratio]{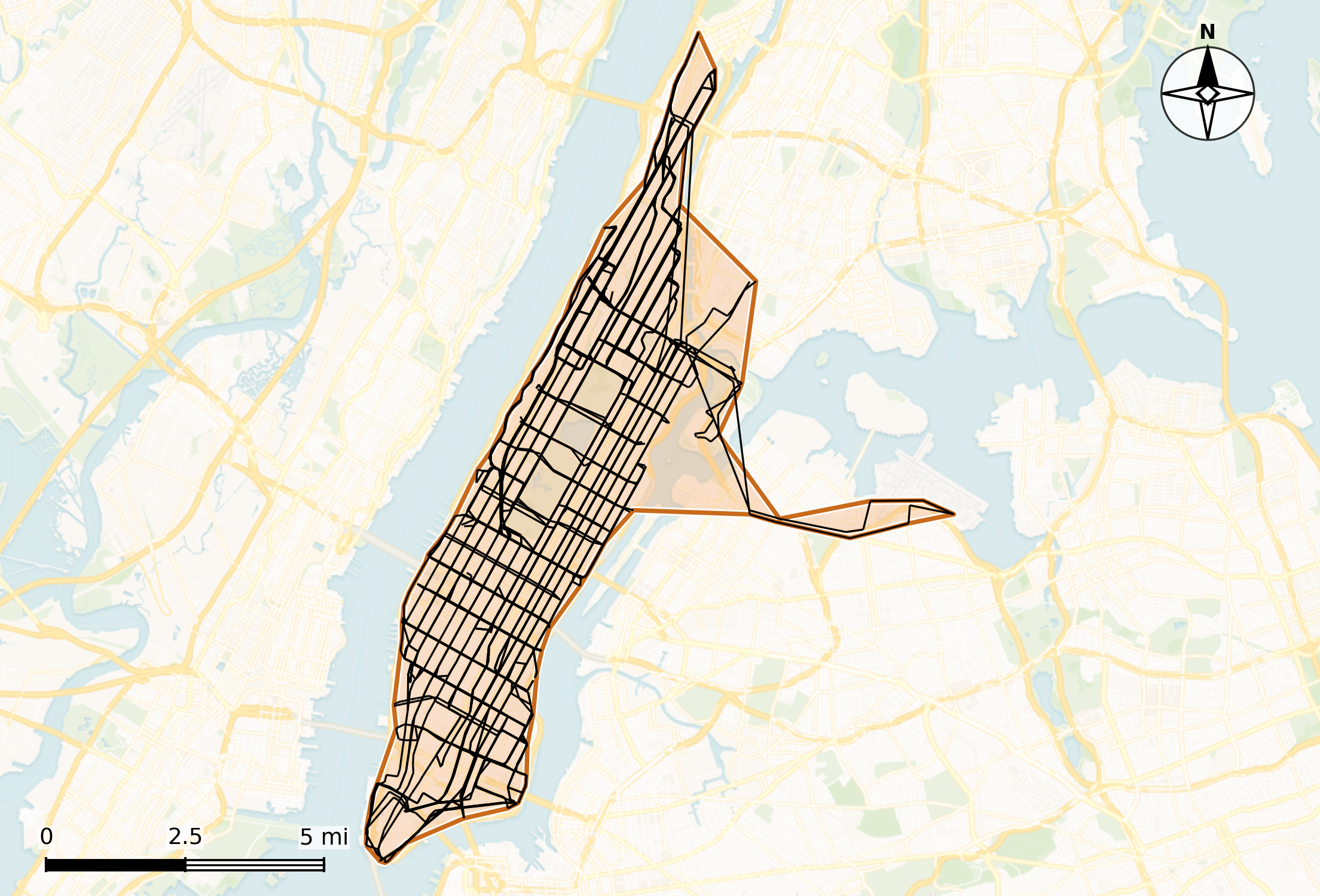}
    \caption{MTA}
    \label{fig:mta_density_map}
\end{subfigure}

\caption{Transit network density maps by agency.}
\label{fig:network_density_maps}
\end{figure}

\Cref{fig:network_density_maps} shows the stop-based network footprints and modeled transit links for the agency instances. \Cref{fig:network_density_bars} summarizes the corresponding transit-link densities.

% \FloatBarrier

% Read Python-generated CSV directly
\pgfplotstableread[col sep=comma]{density/network_density_summary.csv}\densitydata

% Assumes CSV row order from Python is:
% 0 = Santa Monica / BBB
% 1 = Knoxville / KAT
% 2 = Atlanta / MARTA
% 3 = Manhattan / MTA
% 4 = Chicago-CTA / CTA
% 5 = Chicago-PACE / PACE (ignored here)

\pgfplotstablegetelem{0}{Network Density (mi / sq mi)}\of{\densitydata}
\edef\BBBDensity{\pgfplotsretval}

\pgfplotstablegetelem{1}{Network Density (mi / sq mi)}\of{\densitydata}
\edef\KATDensity{\pgfplotsretval}

\pgfplotstablegetelem{2}{Network Density (mi / sq mi)}\of{\densitydata}
\edef\MARTADensity{\pgfplotsretval}

\pgfplotstablegetelem{3}{Network Density (mi / sq mi)}\of{\densitydata}
\edef\MTADensity{\pgfplotsretval}

\pgfplotstablegetelem{4}{Network Density (mi / sq mi)}\of{\densitydata}
\edef\CTADensity{\pgfplotsretval}

\begin{figure}[!ht]
\centering

\begin{tikzpicture}
\begin{axis}[
    width=\textwidth,
    height=0.55\textwidth,
    xbar,
    bar width=12pt,
    xmin=0,
    xmax=40,
    xlabel={Network density (mi/sq mi)},
    xlabel style={font=\small\bfseries},
    tick label style={font=\footnotesize},
    yticklabel style={font=\footnotesize},
    symbolic y coords={BBB,CTA,KAT,MARTA,MTA},
    ytick={BBB,CTA,KAT,MARTA,MTA},
    y dir=reverse,
    enlarge y limits=0.12,
    scaled ticks=false,
    grid=major,
    major grid style={black!10},
    axis line style={black},
    tick align=outside,
    nodes near coords={
        \pgfmathprintnumber[fixed,precision=1]{\pgfplotspointmeta}
    },
    point meta=x,
    every node near coord/.append style={
        font=\footnotesize\bfseries,
        anchor=west,
        xshift=4pt
    },
]

% BBB
\addplot[
    fill=mplblue,
    draw=mplblue!70!black,
    forget plot
]
coordinates {(\BBBDensity,BBB)};

% CTA
\addplot[
    fill=mplorange,
    draw=mplorange!70!black,
    forget plot
]
coordinates {(\CTADensity,CTA)};

% KAT
\addplot[
    fill=mplred,
    draw=mplred!70!black,
    forget plot
]
coordinates {(\KATDensity,KAT)};

% MARTA
\addplot[
    fill=mplpurple,
    draw=mplpurple!70!black,
    forget plot
]
coordinates {(\MARTADensity,MARTA)};

% MTA
\addplot[
    fill=mplbrown,
    draw=mplbrown!70!black,
    forget plot
]
coordinates {(\MTADensity,MTA)};

\end{axis}
\end{tikzpicture}

\caption{Transit network density by agency.}
\label{fig:network_density_bars}
\end{figure}
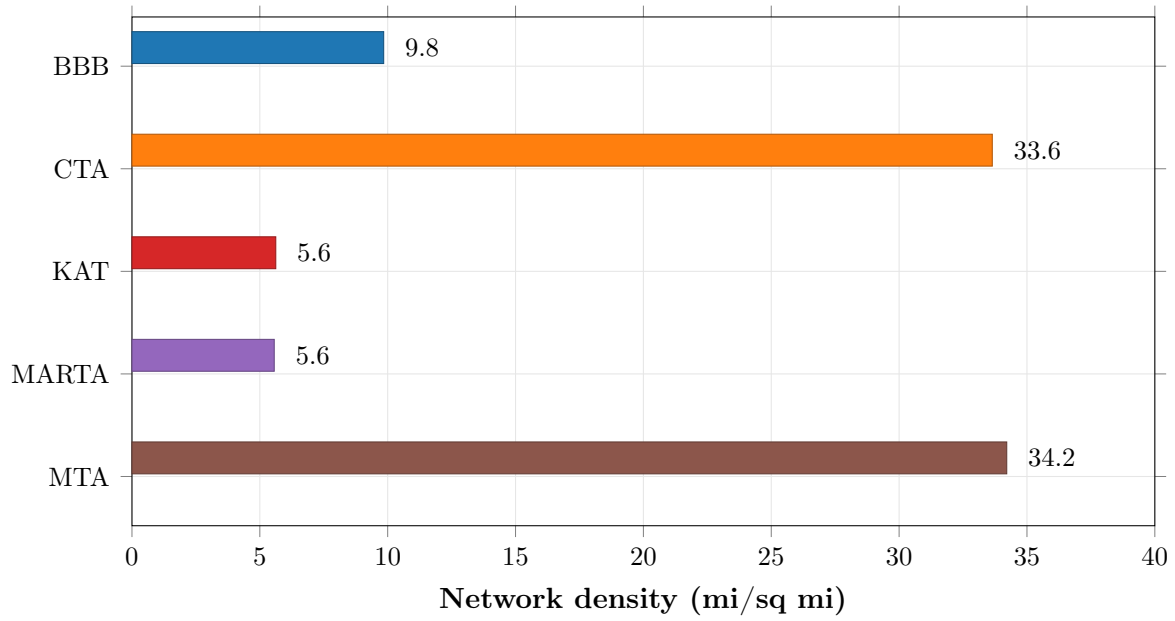

% \FloatBarrier

The transit-link-density values vary substantially across the five analyzed systems. MTA
has the highest transit-link density at approximately 34.2 mi/sq mi, followed closely by
CTA at approximately 33.6 mi/sq mi. BBB has an intermediate density of
approximately 9.8 mi/sq mi, while KAT and MARTA each have densities of
approximately 5.6 mi/sq mi.

At the agency level, the two systems with the highest transit-link densities, CTA and MTA, also have the two lowest average replacement ratios, while the two systems with the lowest densities, KAT and MARTA, have the two highest average ratios. BBB occupies an intermediate position for both measures. This five-agency pattern provides descriptive support for the proposed relationship, although the sample remains too small for statistical inference.

Transit-link density should not be interpreted as a direct measure of BEB schedule feasibility. A possible explanation is that spatially concentrated networks may offer shorter or more compatible connections among trips, garages, and charging locations \cite{adler_mirchandani_2017,li_2014,liu_ceder_2020,bazarnovi_2026}. However, the density measure does not directly quantify those movements. Charger availability, timetable structure, layover distribution, and garage assignment may also influence the optimized replacement ratio \cite{alvo_2021,bao_2023,duan_2023,bazarnovi_2026}. The vehicle activity results provide schedule-level evidence that can be considered alongside the spatial comparison.

\subsection{Vehicle Activity Time Allocation}

\begin{figure}[!ht]
\centering

\begin{subfigure}{0.48\textwidth}
    \centering
    \includegraphics[width=\textwidth]{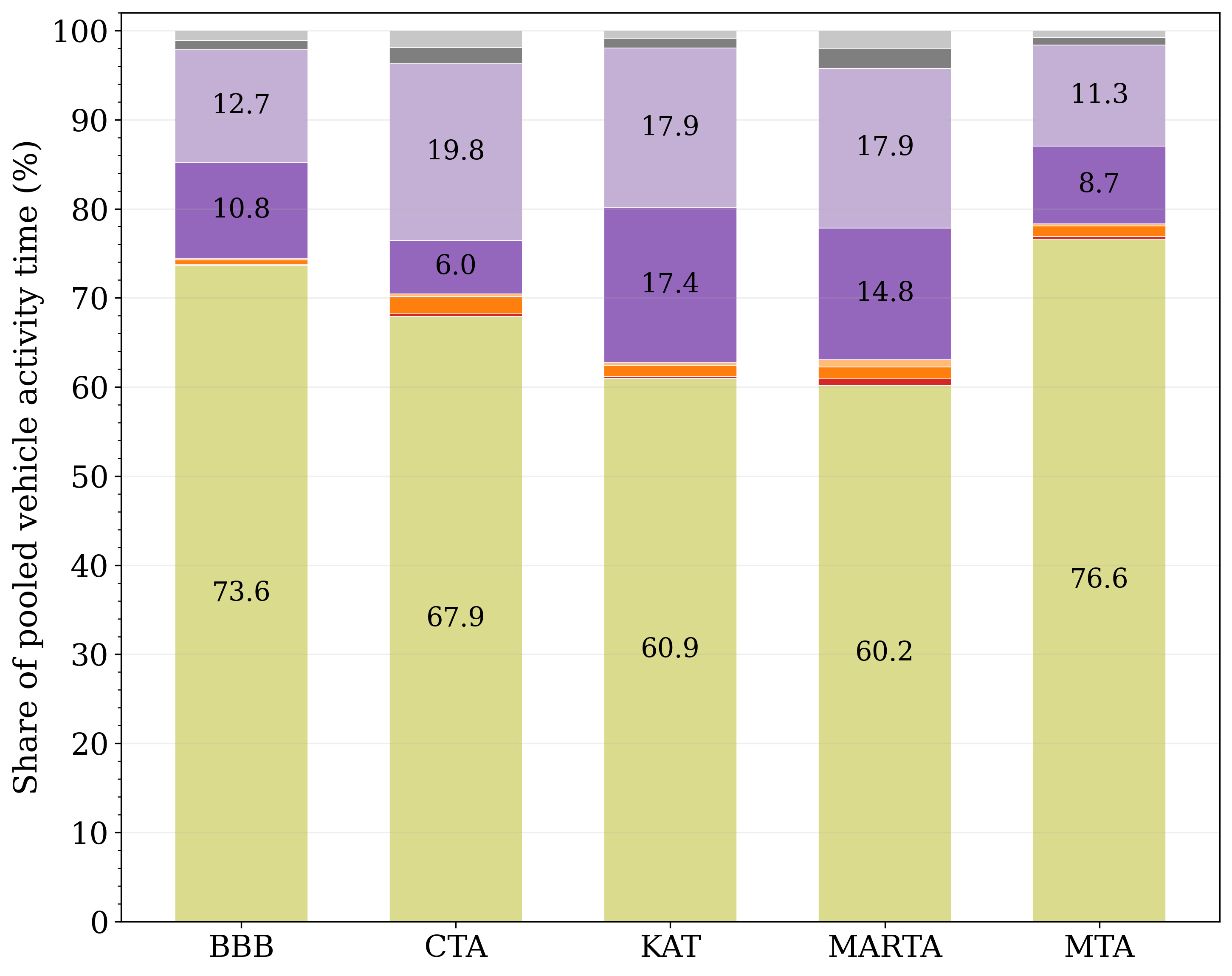}
    \caption{Diesel buses}
    \label{fig:activity_allocation_db}
\end{subfigure}
\hfill
\begin{subfigure}{0.48\textwidth}
    \centering
    \includegraphics[width=\textwidth]{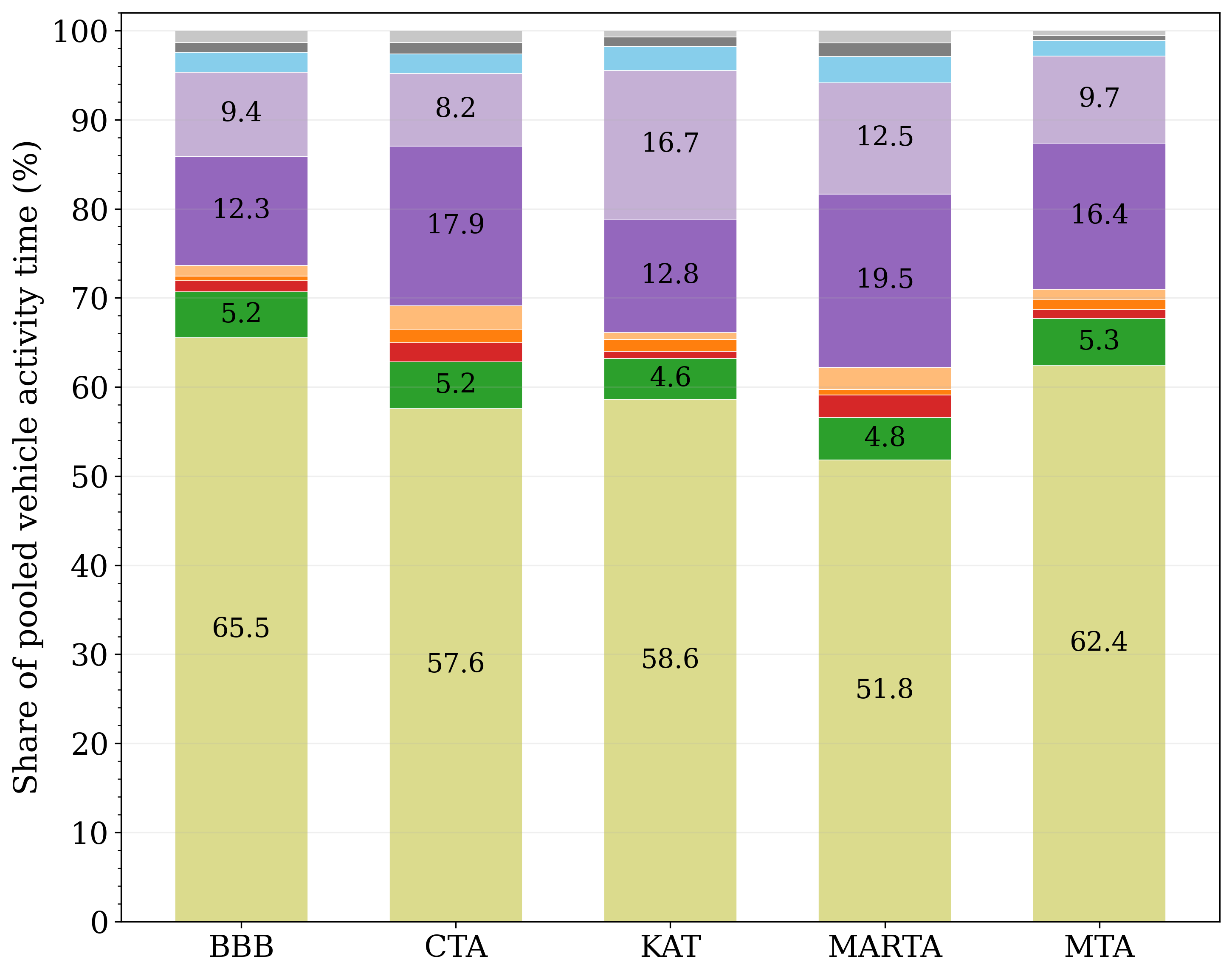}
    \caption{Battery electric buses}
    \label{fig:activity_allocation_beb}
\end{subfigure}

\vspace{0.5em}

\includegraphics[width=0.95\textwidth]{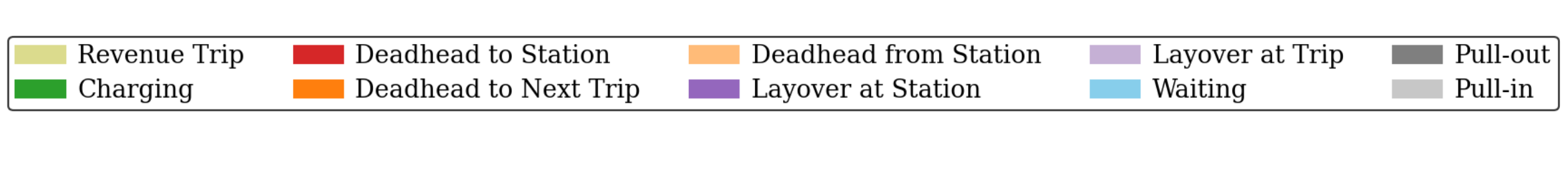}

\caption{Pooled vehicle activity-time allocation for DBs and BEBs across mixed-fleet scenarios.}
\label{fig:activity_allocation}
\end{figure}

\Cref{fig:activity_allocation} shows the pooled activity-time allocation for DBs and BEBs across all mixed-fleet scenarios containing BEBs. The DB panel represents DBs retained in the optimized mixed fleets rather than vehicles from the diesel-only baseline.

% \FloatBarrier

The DB activity allocations differ across the five analyzed agencies. BBB and MTA have the highest DB revenue-time shares, at approximately 74\% and 77\%, respectively, while CTA occupies an intermediate position with approximately 68\% of pooled retained-DB scheduled activity time assigned to revenue trips. KAT and MARTA have lower DB revenue-time shares, at approximately 61\% and 60\%, respectively, along with greater allocations to layover and other non-revenue activities.

Differences in BEB activity allocation may help explain why some agencies require more BEBs per DB displaced. CTA allocates approximately 58\% of pooled scheduled activity time to revenue trips and approximately 5\% to charging. BBB has the highest BEB revenue-time share among the five agencies, at approximately 66\%, followed by MTA at approximately 62\%. KAT and MARTA have lower BEB revenue-time shares of approximately 59\% and 52\%, respectively. These patterns broadly align with the replacement-ratio results: CTA and MTA have lower average ratios, KAT and MARTA have higher ratios, and BBB occupies an intermediate position.

Replacement efficiency is not determined by charging time alone. Charging accounts for a relatively narrow range of approximately 4.6\% to 5.3\% of pooled BEB scheduled activity time across the five agencies, whereas larger differences are observed in revenue and layover shares. Charging, deadheading, waiting, pull-out, and pull-in activities reduce the portion of a vehicle's schedule available for revenue service \cite{cokyasar_2023_large_scale,davatgari_2024_ndo,bazarnovi_2026}. Layover also reduces the revenue share, although it may provide usable time for charging and should not automatically be treated as an avoidable burden \cite{hu_2021,duan_2023,bazarnovi_2026}.

Taken together, the spatial and activity results provide possible explanations for the observed replacement-ratio differences. CTA and MTA combine comparatively concentrated networks with lower replacement ratios, while MARTA and KAT have lower transit-link densities, lower BEB revenue-time shares, and higher ratios. BBB has a comparatively high BEB revenue-time share but an intermediate replacement ratio, indicating that revenue-time share alone does not fully explain replacement efficiency. These relationships are descriptive, and the five-agency sample does not establish that transit-link density or any single activity category causes the observed outcomes.

% =============================================================================
\section{Conclusion}\label{sec:conclusion}
% =============================================================================

This study evaluated BEB replacement efficiency using outputs from the integrated mixed BEB/DB scheduling and flexible-charging model developed by \citet{bazarnovi_2026}. Five transit systems---BBB, CTA, KAT, MARTA, and MTA---were evaluated under BEB/DB purchase-cost ratios from 1.0 to 1.4.

The optimized fleets remained majority-electric, but replacement was not one-to-one. Average replacement ratios were 1.101 for CTA, 1.171 for MTA, 1.201 for BBB, 1.234 for MARTA, and 1.245 for KAT. These findings are consistent with prior electric-vehicle scheduling studies showing that operationally feasible electrification may require additional vehicles because of range, charging, and non-revenue activity constraints \cite{cokyasar_2023_large_scale,davatgari_2024_ndo,liu_ceder_2020,bazarnovi_2026}.

Electric fleet share and replacement ratio describe different dimensions of the optimized outcome. Fleet share describes composition, whereas the replacement ratio indicates how many BEBs are selected per DB displaced. CTA and MTA combine higher transit-link densities with lower replacement ratios, while MARTA and KAT have lower densities, lower revenue-time shares, and higher ratios. This pattern suggests that network structure and activity allocation may influence replacement efficiency, but the relationships are descriptive rather than causal.

Agencies should therefore evaluate replacement efficiency, total fleet requirements, and vehicle activity allocation alongside electric fleet share rather than assume one-to-one substitution. The analysis is limited to five systems, fixed timetables and garage assignments, garage-only charging and service-based charger-capacity assumptions, and a visually calibrated stop-based network footprint. Future work should examine additional agencies and directly measure inter-trip deadheading, garage and charger accessibility, charger capacity, garage space, and infrastructure investment \cite{davatgari_2024,bazarnovi_2025_charging_location}.

\section*{Acknowledgments}
This material is based on work supported by the U.S. Department of Energy (DOE), Office of Science, under contract number DE-AC02-06CH11357. This report and the work described were sponsored by the DOE Transportation Technologies Office (TTO) under the Pathways to Affordable, Convenient, and Efficient Regional Mobility, an initiative of the Energy Efficient Mobility Systems (EEMS) Program. Erin Boyd, a DOE Office of Critical Minerals and Energy Innovation (CMEI) manager, played an important role in establishing the project concept, advancing implementation, and providing guidance.

\section*{Author Contributions} 
The authors confirm contribution to the paper as follows: study conception and design: Bhidya, H., Cokyasar, T., Verbas, O.; data collection: Bhidya, H., Cokyasar, T.; methodology: Bhidya, H., Cokyasar, T., Verbas, O.; software: Bhidya, H., Cokyasar, T.; validation: Bhidya, H., Cokyasar, T., Verbas, O.; analysis and interpretation of results: Bhidya, H., Cokyasar, T., Verbas, O.; visualization: Bhidya, H.; draft manuscript preparation: Bhidya, H. All authors reviewed the results and approved the final version of the manuscript.

\section*{AI Use} 
During the preparation of this manuscript, the authors used Claude Opus 4.6 and Claude Sonnet 4.6 through the ARGO AI platform at Argonne National Laboratory to improve language, spelling, readability, literature review, coding, and result analysis. After using these tools, the authors meticulously reviewed and edited the entire text, code, and references (including manuscripts and all online links) and take full responsibility for the content of the publication.

\section*{Declaration of Conflict of Interest} 
None declared.

\section*{ORCID}
\def\orcid#1{\kern .08em\href{https://orcid.org/#1}{\includegraphics[keepaspectratio,width=0.7em]{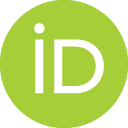}}}
\noindent Hadi Bhidya \orcid{0009-0005-5483-4603} \url{https://orcid.org/0009-0005-5483-4603}\\
Taner Cokyasar \orcid{0000-0001-9687-6725} \url{https://orcid.org/0000-0001-9687-6725}\\
Omer Verbas \orcid{0000-0002-4831-7475} \url{https://orcid.org/0000-0002-4831-7475}\\

\newpage
\bibliographystyle{elsarticle-harv} 
\bibliography{trb_references}

\end{document}